\documentclass[11pt]{amsart}
\usepackage{amsfonts}
\usepackage{amsmath,amssymb,amsthm}
\usepackage{color,graphics}
\usepackage{srcltx}
\usepackage{bm}
\usepackage{amsmath}
\usepackage{amssymb}
\usepackage{graphicx}
\usepackage{amsmath,amsfonts,amssymb,amsthm,graphics,amscd,enumerate,verbatim}
\usepackage{color}
\usepackage{url}

\newtheorem{theorem}{Theorem}[section]
\newtheorem{corollary}[theorem]{Corollary}
\newtheorem{lemma}[theorem]{Lemma}
\newtheorem{proposition}[theorem]{Proposition}
\newtheorem{definition}{Definition}

\newtheorem{remark}[theorem]{Remark}

\input{tcilatex}

\begin{document}
\title{On some generalized inverses and partial orders in $\ast $-rings}
\author{Janko Marovt}
\address[Janko Marovt]{Janko Marovt, University of Maribor, Faculty of
Economics and Business, Razlagova 14, SI-2000 Maribor, Slovenia, and IMFM,
Jadranska 19, SI-1000 Ljubljana, Slovenia}
\email{janko.marovt@um.si}
\author{Dijana Mosi\'c}
\address[Dijana Mosi\'c]{Faculty of Sciences and Mathematics, University of
Ni\v s, P.O. Box 224, 18000 Ni\v s, Serbia}
\email{dijana@pmf.ni.ac.rs}
\author{Insa Cremer}
\address[Insa Cremer]{Department of Mathematics, University of Latvia,
Jelgavas 3, Riga, 1004, Latvia}
\email{insa.kremere@gmail.com}
\thanks{This first author acknowledges the financial support from the
Slovene Research Agency, ARRS (research core funding No. P1-0288). The second
author is supported by the Ministry of Education, Science and Technological
Development, Republic of Serbia (grant no. 174007
(451-03-9/2021-14/200124)). The authors also acknowledge the bilateral
projects between Serbia and Slovenia (Generalized inverses, operator
equations and applications), and Slovenia and Latvia (Ordered structures in
Rickart rings) were financially supported by the Ministry of Education,
Science and Technological Development, Republic of Serbia (grant no.
337-00-21/2020-09/32) and by the Slovene Research Agency, ARRS (grants
BI-RS/20-21-039 and BI-LV/20-22-002).}

\begin{abstract}
Let $\mathcal{R}$ be a unital ring with involution. The notions of
1MP-inverse and MP1-inverse are extended from $M_{m,n}(\mathbb{C)}$, the set
of all $m\times n $ matrices over $\mathbb{C}$, to the set $\mathcal{R}%
^{\dagger}$ of all Moore-Penrose invertible elements in $\mathcal{R}$. We
study partial orders on $\mathcal{R}^{\dagger}$ that are induced by
1MP-inverses and MP1-inverses. We also extend to the setting of Rickart $\ast
$-rings the concept of another partial order, called the plus order, which
has been recently introduced on the set of all bounded linear operators
between Hilbert spaces. Properties of these relations are investigated and
some known results are thus generalized.
\end{abstract}

\maketitle

\section{Introduction}

Let $\mathcal{R}$ be a $\ast $-ring, i.e. a ring equipped with involution $%
\ast $. We call an element $a\in \mathcal{R}$ a $\ast $\textit{-regular} or
\textit{Moore-Penrose invertible }with respect to $\ast $ if there exists $%
x\in \mathcal{R}$ that satisfies the following four equations:
\begin{equation}
axa=a,\quad xax=x,\quad (ax)^{\ast }=ax,\quad (xa)^{\ast }=xa.
\label{def_MP_inverse}
\end{equation}%
If such $x$ exists, we write $x=a^{\dagger }$ and call it the\textit{\
Moore-Penrose inverse} of $a$. It is known that $a^{\dagger }$ is unique if
it exists. We denote the set of all $\ast $-regular elements in $\mathcal{R}$
by $\mathcal{R}^{\dagger }$. A ring $\mathcal{R}$ where every element is $%
\ast $-regular is called a $\ast $\textit{-regular ring}. An example of a $%
\ast $-regular ring is the set $M_{n}(\mathbb{C})$ of all complex $n\times n$
matrices where $A^{\ast }$ denotes the conjugate transpose of $A\in M_{n}(%
\mathbb{C})$. There are many other generalized inverses that can be
introduced on $\mathcal{R}$ and some of them may be defined with a subset of
the set of equations (\ref{def_MP_inverse}). We say that an element $a\in
\mathcal{R}$ is \textit{regular} if there exists $x\in \mathcal{R}$ that
satisfies the first equation in (\ref{def_MP_inverse}). Such $x$, if it
exists, is called an\textit{\ inner generalized inverse} or $\{1\}$\textit{%
-inverse} of $a$, and we write $x=a^{-}$, i.e. $aa^{-}a=a$. The set of all
regular elements in $\mathcal{R}$ is denoted by $\mathcal{R}^{(1)}$. A ring $%
\mathcal{R}$ is called \textit{von Neumann regular} if every element in $%
\mathcal{R}$ is regular. If $x\in \mathcal{R}$ satisfies the first three
equations in (\ref{def_MP_inverse}), then it is called a $\{1,2,3\}$\textit{%
-inverse} of $a$. We denote the sets of all $\{1\}$-inverses and $\{1,2,3\}$%
-inverses of $a\in \mathcal{R}$ by $a\{1\}$ and $a\{1,2,3\}$, respectively.

In a recent paper \cite{Hernandez}, Hern\'{a}ndez et al. introduced two new
types of hybrid generalized inverses on the set $M_{m,n}(\mathbb{C})$ of all
complex $m\times n$ matrices. Note here that for $A\in M_{m,n}(\mathbb{C})$
the definitions of its inner generalized inverse $A^{-}\in M_{n,m}(\mathbb{C)}$ and
its unique Moore-Penrose inverse $A^{\dag }\in M_{n,m}(\mathbb{C})$ are the
same as on $\mathcal{R}$.

\begin{definition}
Let $A\in M_{m,n}(\mathbb{C})$. For each inner generalized inverse $A^{-}$
of $A$, the matrices
\begin{equation*}
A^{-\dagger }=A^{-}AA^{\dagger }\text{\quad and\quad }A^{\dagger
-}=A^{\dagger }AA^{-}
\end{equation*}%
are called a 1MP-inverse and a MP1-inverse of $A$, respectively.
\end{definition}

An inner generalized inverse and the Moore-Penrose inverse of $A\in M_{m,n}(%
\mathbb{C)}$ always exist and this guarantees the existence of a 1MP- and a
MP1-inverse for every $A\in M_{m,n}(\mathbb{C)}$. In general, 1MP- and
MP1-inverses are not unique.

In \cite{Hernandez}, the 1MP- and MP1-inverses of $A\in M_{m,n}(\mathbb{C)}$
were characterized. Also, a partial order induced by a 1MP-inverse and its
dual case associated with a MP1-inverse were introduced. It is the aim of
this paper to generalize the concept of a 1MP-inverse and a MP1-inverse to $%
\ast $-regular rings or more generally to the set of all Moore-Penrose
invertible elements in a $\ast $-ring $\mathcal{R}$, study their properties,
and the properties of partial orders on $\mathcal{R}$ associated with these
generalized inverses. We thus extend the results of \cite{Hernandez}. The
paper is structured as follows. In Section 2, we give some preliminaries and
in Section 3 we introduce the concept of 1MP-inverses in the setting of $%
\ast $-rings and present some characterizations of these generalized
inverses. We introduce a relation associated with 1MP-inverses and study its
properties in Section 4. The dual case (the MP1-inverses and a relation
associated with these generalized inverses) is studied in Section 5. In the
last section, we extend to the setting of Rickart $\ast$-rings the concept
of another partial order, called the plus order, which has been recently
introduced in \cite{Arias} on the set $B(\mathcal{H},\mathcal{K})$ of all
bounded linear operators between Hilbert spaces $\mathcal{H}$ and $\mathcal{K%
}$.

\section{Preliminaries}

Unless stated otherwise, let from now $\mathcal{R}$ be a $\ast $-ring
with the (multiplicative) identity $1$. If for $p\in \mathcal{R}$, $p^{2}=p$%
, then $p$ is said to be an idempotent. A projection $p\in \mathcal{R}$ is a
self-adjoint idempotent, i.e. $p=p^{2}=p^{\ast }$. The equality $%
1=e_{1}+e_{2}+\cdots +e_{n}$, where $e_{1},e_{2},\dots ,e_{n}$ are
idempotents in $\mathcal{R}$ and $e_{i}e_{j}=0$ for $i\neq j$, is called a
decomposition of the identity of $\mathcal{R}$. Let $1=e_{1}+e_{2}+\dots
+e_{n}$ and $1=f_{1}+f_{2}+\dots +f_{n}$ be two decompositions of the
identity of $\mathcal{R}$. We have
\begin{equation*}
x=1\cdot x\cdot 1=(e_{1}+e_{2}+\dots +e_{n})x(f_{1}+f_{2}+\dots
+f_{n})=\sum_{i,j=1}^{n}e_{i}xf_{j}.
\end{equation*}%
Then any $x\in \mathcal{R}$ can be uniquely represented in the following
matrix form:
\begin{equation}
x=\left[
\begin{array}{ccc}
x_{11} & \cdots & x_{1n} \\
\vdots & \ddots & \vdots \\
x_{n1} & \cdots & x_{nn}%
\end{array}%
\right] _{e\times f},  \label{Matrix_formulation}
\end{equation}%
where $x_{ij}=e_{i}xf_{j}\in e_{i}\mathcal{R}f_{j}$. With $e\times f$ we
emphasize the use of the decompositions of the identity $1=e_{1}+e_{2}+%
\cdots +e_{n}$ on the left side and $1=f_{1}+f_{2}+\dots +f_{n}$ on the
right side of $x=1\cdot x\cdot 1$. If $x=(x_{ij})_{e\times f}$ and $%
y=(y_{ij})_{e\times f}$, then $x+y=(x_{ij}+y_{ij})_{e\times f}$. Moreover,
if $1=g_{1}+\cdots +g_{n}$ is another decomposition of the identity of $%
\mathcal{R}$ and $z=(z_{ij})_{f\times g}$, then, by the orthogonality of the
idempotents involved, $xz=\left( \sum_{k=1}^{n}x_{ik}z_{kj}\right) _{e\times
g}$. Thus, if we have decompositions of the identity of $\mathcal{R}$, then
the usual algebraic operations in $\mathcal{R}$ can be interpreted as simple
operations between appropriate $n\times n$ matrices over $\mathcal{R}$. When
$n=2$ and $p,q\in \mathcal{R}$ are idempotents, we may write
\begin{equation*}
x=pxq+px(1-q)+(1-p)xq+(1-p)x(1-q)=\left[
\begin{array}{cc}
x_{11} & x_{12} \\
x_{21} & x_{22}%
\end{array}%
\right] _{p\times q}.
\end{equation*}%
Here $x_{11}=pxq,$ $x_{12}=px(1-q),$ $x_{21}=(1-p)xq,$ $x_{22}=(1-p)x(1-q)$.

By (\ref{Matrix_formulation}) we may write
\begin{equation*}
x^{\ast }=\left[
\begin{array}{ccc}
x_{11}^{\ast } & \cdots & x_{n1}^{\ast } \\
\vdots & \ddots & \vdots \\
x_{1n}^{\ast } & \cdots & x_{nn}^{\ast }%
\end{array}%
\right] _{f^{\ast }\times e^{\ast }},
\end{equation*}%
where this matrix representation of $x^{\ast }$ is given relative to the
decompositions of the identity $1=f_{1}^{\ast }$ $+\dots +f_{n}^{\ast }$ and
$1=e_{1}^{\ast }+\dots +e_{n}^{\ast }$ .

Let $a\in \mathcal{R}$ and $a^{\circ }$ denote the right annihilator of $a$,
i.e. the set $a^{\circ }=\{x\in \mathcal{R}:ax=0\}$. Similarly we denote the
left annihilator $^{\circ }a$ of $a,$ i.e. the set $^{\circ }a=\{x\in
\mathcal{R}:xa=0\}$. Suppose that $p,q\in \mathcal{R}$ are such idempotents
that $^{\circ }a=$ $^{\circ }p$ and $a^{\circ }=q^{\circ }$. Observe (or see
\cite[Lemma 2.2]{DjordjevicRakicMarovt}) that $^{\circ }p=\mathcal{R(}1-p%
\mathcal{)}$ and $q^{\circ }=(1-q)\mathcal{R}$. It follows that then $a=paq$%
, i.e.%
\begin{equation}
a=\left[
\begin{array}{cc}
a & 0 \\
0 & 0%
\end{array}%
\right] _{p\times q}.  \label{Matrix form_basic}
\end{equation}%
Suppose $a\in \mathcal{R}^{(1)}$, fix $h\in a\{1\}$, and let $p=ah$ and $q=ha
$. Clearly, $p$ and $q$ are idempotents. Moreover, $pa=a$ and $aq=a$, and so
$^{\circ }a=$ $^{\circ }p$ and $a^{\circ }=q^{\circ }$. We may thus write $a$
in the matrix form (\ref{Matrix form_basic}). Also, by \cite[page 1044]%
{RakicDecomposition}, $k\in a\{1\}$ if and only if%
\begin{equation}
k=\left[
\begin{array}{cc}
hah & k_{12} \\
k_{21} & k_{22}%
\end{array}%
\right] _{q\times p}  \label{matrix_from_generalized_inner}
\end{equation}%
where $k_{12}\in q\mathcal{R}(1-p)$, $k_{21}\in (1-q)\mathcal{R}p$, and $%
k_{22}\in (1-q)\mathcal{R}(1-p)$.

Suppose that $a\in \mathcal{R}^{\dagger }$
and let now $p=aa^{\dagger }$ and $q=a^{\dagger }a$. Then, $p$ and $q$ are
projections and $^{\circ }a=$ $^{\circ }p$ and $a^{\circ }=q^{\circ }$.
Also, since $a^{\dagger }p=a^{\dagger }=$ $qa^{\dagger }$, we may conclude
that $^{\circ }a^{\dagger }=$ $^{\circ }q$ and $a^{\dagger \circ }=p^{\circ }
$. So,
\begin{equation*}
a^{\dagger }=\left[
\begin{array}{cc}
a^{\dagger } & 0 \\
0 & 0%
\end{array}%
\right] _{q\times p}.
\end{equation*}

\section{1MP-inverses in a $\ast $-ring}

Let us now extend the concept of 1MP-inverses to the set of all
Moore-Penrose invertible elements in a $\ast $-ring $\mathcal{R}$.

\begin{definition}
Let $a\in \mathcal{R}^{\dagger }$. For each $a^{-}\in $ $a\{1\}$, the
element
\begin{equation*}
a^{-\dagger }=a^{-}aa^{\dagger }
\end{equation*}%
is called a 1MP-inverse of $a$. The set of all 1MP-inverses of $a$ is
denoted by $a\{-\dagger \}$.
\end{definition}

Since $a^{\dagger }=a^{\dagger }aa^{\dagger }$, we observe that $a^{\dagger
}\in a\{-\dagger \}$. So, for every $a\in \mathcal{R}^{\dagger }$, the set $%
a\{-\dagger \}$ is nonempty. Let $a\in \mathcal{R}^{\dagger }$ and $a^{-}\in
$ $a\{1\}$. Denote $p=aa^{\dagger }$ and $q=a^{\dagger }a$. Then
\begin{equation*}
a=\left[
\begin{array}{cc}
a & 0 \\
0 & 0%
\end{array}%
\right] _{p\times q}\quad \text{and\quad }a^{\dagger }=\left[
\begin{array}{cc}
a^{\dagger } & 0 \\
0 & 0%
\end{array}%
\right] _{q\times p}.
\end{equation*}%
Observe that by (\ref{matrix_from_generalized_inner})
\begin{equation*}
a^{-}=\left[
\begin{array}{cc}
a^{\dagger }aa^{\dagger } & d_{12} \\
d_{21} & d_{22}%
\end{array}%
\right] _{q\times p}=\left[
\begin{array}{cc}
a^{\dagger } & d_{12} \\
d_{21} & d_{22}%
\end{array}%
\right] _{q\times p}
\end{equation*}%
for some $d_{12}\in q\mathcal{R}(1-p)$, $d_{21}\in (1-q)\mathcal{R}p$, and $%
d_{22}\in (1-q)\mathcal{R}(1-p)$. Thus,
\begin{eqnarray*}
a^{-}aa^{\dagger } &=&\left[
\begin{array}{cc}
a^{\dagger } & d_{12} \\
d_{21} & d_{22}%
\end{array}%
\right] _{q\times p}\left[
\begin{array}{cc}
a & 0 \\
0 & 0%
\end{array}%
\right] _{p\times q}\left[
\begin{array}{cc}
a^{\dagger } & 0 \\
0 & 0%
\end{array}%
\right] _{q\times p} \\
&=&\left[
\begin{array}{cc}
a^{\dagger } & d_{12} \\
d_{21} & d_{22}%
\end{array}%
\right] _{q\times p}\left[
\begin{array}{cc}
p & 0 \\
0 & 0%
\end{array}%
\right] _{p\times p} \\
&=&\left[
\begin{array}{cc}
a^{\dagger }p & 0 \\
d_{21}p & 0%
\end{array}%
\right] _{q\times p}=\left[
\begin{array}{cc}
a^{\dagger } & 0 \\
d_{21} & 0%
\end{array}%
\right] _{q\times p}.
\end{eqnarray*}%
Since $\left[
\begin{array}{cc}
a^{\dagger } & d_{12} \\
d_{21} & d_{22}%
\end{array}%
\right] _{q\times p}\in a\{1\}$ for any $d_{12}\in q\mathcal{R}(1-p)$, $%
d_{21}\in (1-q)\mathcal{R}p$, and $d_{22}\in (1-q)\mathcal{R}(1-p)$, we may
conclude that
\begin{equation}
a\{-\dagger \}=\left\{ \left[
\begin{array}{cc}
a^{\dagger } & 0 \\
d_{21} & 0%
\end{array}%
\right] _{q\times p}:d_{21}\in (1-q)\mathcal{R}p\text{ is arbitrary}\right\}
.  \label{minus plus set}
\end{equation}

We will next give some characterizations of 1MP-inverses. First, let us
prove a simple auxiliary result.

\begin{lemma}
\label{Lemma_first}Let $a\in \mathcal{R}^{\dagger }$. Then $a\{-\dagger
\}\subseteq a\{1,2,3\}$, $a^{-\dagger }a=a^{-}a$, and $aa^{-\dagger
}=aa^{\dagger }.$
\end{lemma}

\begin{proof}
Let $a^{-\dagger }\in a\{-\dagger \}$, i.e. $a^{-\dagger }=a^{-}aa^{\dagger
} $ for some $a^{-}\in a\{1\}$. From
\begin{eqnarray*}
aa^{-\dagger }a &=&aa^{-}aa^{\dagger }a=a, \\
a^{-\dagger }aa^{-\dagger } &=&a^{-}aa^{\dagger }aa^{-}aa^{\dagger
}=a^{-}aa^{\dagger }=a^{-\dagger },
\end{eqnarray*}%
and
\begin{equation*}
\left( aa^{-\dagger }\right) ^{\ast }=\left( aa^{-}aa^{\dagger }\right)
^{\ast }=\left( aa^{\dagger }\right) ^{\ast }=aa^{\dagger
}=aa^{-}aa^{\dagger }=aa^{-\dagger },
\end{equation*}%
it follows that $a^{-\dagger }\in a\{1,2,3\}$ and thus $a\{-\dagger
\}\subseteq a\{1,2,3\}$. Also, $aa^{-\dagger }=$ $aa^{\dagger }$ and $%
a^{-\dagger }a=a^{-}aa^{\dagger }a=a^{-}a$.
\end{proof}

\begin{remark}
Let $a\in \mathcal{R}^{\dagger }$. Since $aa^{\dagger }$ is a projection and
$a^{-}a$ is an idempotent, we may conclude, by Lemma \ref{Lemma_first},
that $aa^{-\dagger }$ is a projection and $a^{-\dagger }a$ is an idempotent.
Moreover, from $a=\left( aa^{\dagger }\right) a$ and since for $z\in
\mathcal{R}$, $za=0$ implies $zaa^{\dagger }=0$, we have that
$^{\circ }\left( aa^{-\dagger }\right) =$ $^{\circ }\left( aa^{\dagger
}\right) =$ $^{\circ }a$. Similarly we obtain that $\left( a^{-\dagger
}a\right) ^{\circ }=a^{\circ }$.
\end{remark}

By Lemma \ref{Lemma_first} it follows that every $x=a^{-\dagger }\in
a\{-\dagger \}$ satisfies the following system of equations:
\begin{equation}
xax=x\quad \text{and\quad }ax=aa^{\dagger }.  \label{System}
\end{equation}

With the next theorem we extend \cite[Theorem 3.1]{Hernandez} from the set
of all $m\times n$ complex matrices to the set $\mathcal{R}^{\dagger }$ in a
$\ast $-ring $\mathcal{R}$ with identity.

\begin{theorem}
\label{Theorem1}Let $a\in \mathcal{R}^{\dagger }$. Then the following
statements are equivalent:

\begin{itemize}
\item[(i)] $z\in a\{-\dagger \}$;

\item[(ii)] $z\in \mathcal{R}$ is a solution of the system $(\ref{System})$;

\item[(iii)] $z\in a\{1,2,3\}$.
\end{itemize}
\end{theorem}

\begin{proof}
(i)$\Rightarrow $(ii): It follows by Lemma \ref{Lemma_first}.

(ii)$\Rightarrow $(i): Assume that for $z\in \mathcal{R}$ we have $zaz=z$
and $az=aa^{\dagger }$. Let $p=aa^{\dagger }$ and $q=a^{\dagger }a$. Then
\begin{equation*}
a=\left[
\begin{array}{cc}
a & 0 \\
0 & 0%
\end{array}%
\right] _{p\times q}\quad \text{and\quad }a^{\dagger }=\left[
\begin{array}{cc}
a^{\dagger } & 0 \\
0 & 0%
\end{array}%
\right] _{q\times p}.
\end{equation*}%
Let $z=\left[
\begin{array}{cc}
z_{11} & z_{12} \\
z_{21} & z_{22}%
\end{array}%
\right] _{q\times p}$. Then
\begin{equation*}
az=\left[
\begin{array}{cc}
a & 0 \\
0 & 0%
\end{array}%
\right] _{p\times q}\left[
\begin{array}{cc}
z_{11} & z_{12} \\
z_{21} & z_{22}%
\end{array}%
\right] _{q\times p}=\left[
\begin{array}{cc}
az_{11} & az_{12} \\
0 & 0%
\end{array}%
\right] _{p\times p}
\end{equation*}%
and
\begin{equation*}
aa^{\dagger }=p=\left[
\begin{array}{cc}
p & 0 \\
0 & 0%
\end{array}%
\right] _{p\times p}.
\end{equation*}%
So, $az_{12}=0$ and since $z_{12}\in q\mathcal{R}(1-p)$ and $a^{\circ
}=q^{\circ }$ this implies that $z_{12}=qz_{12}=0$. Also, $az_{11}=p$, i.e. $%
az_{11}=aa^{\dagger }$. Since $z_{11}\in q\mathcal{R}p$ we therefore get
\begin{equation*}
z_{11}=qz_{11}=a^{\dagger }az_{11}=a^{\dagger }aa^{\dagger }=a^{\dagger }.
\end{equation*}%
From $z=zaz$ we obtain
\begin{eqnarray*}
z &=&\left[
\begin{array}{cc}
a^{\dagger } & 0 \\
z_{21} & z_{22}%
\end{array}%
\right] _{q\times p}\left[
\begin{array}{cc}
a & 0 \\
0 & 0%
\end{array}%
\right] _{p\times q}\left[
\begin{array}{cc}
a^{\dagger } & 0 \\
z_{21} & z_{22}%
\end{array}%
\right] _{q\times p} \\
&=&\left[
\begin{array}{cc}
a^{\dagger } & 0 \\
z_{21} & z_{22}%
\end{array}%
\right] _{q\times p}\left[
\begin{array}{cc}
p & 0 \\
0 & 0%
\end{array}%
\right] _{p\times p} \\
&=&\left[
\begin{array}{cc}
a^{\dagger }p & 0 \\
z_{21}p & 0%
\end{array}%
\right] _{q\times p}=\left[
\begin{array}{cc}
a^{\dagger } & 0 \\
z_{21} & 0%
\end{array}%
\right] _{q\times p}.
\end{eqnarray*}%
By (\ref{minus plus set}) we may conclude that $z\in a\{-\dagger \}$.

(ii)$\Rightarrow $(iii): Let $z\in \mathcal{R}$ be a solution of the system $%
(\ref{System})$. Since statement (ii) implies statement (i), $z\in
a\{-\dagger \}$ and thus by Lemma \ref{Lemma_first} $z\in a\{1,2,3\}$.

(iii)$\Rightarrow $(ii): Let $z\in a\{1,2,3\}$, i.e. $aza=a$, $zaz=z$, and $%
\left( az\right) ^{\ast }=az$. Let $p=aa^{\dagger }$ and $q=a^{\dagger }a$.
To conclude the proof we must show that $az=aa^{\dagger }=p$. Since $z\in
a\{1\}$, we may write by (\ref{matrix_from_generalized_inner})
\begin{equation*}
z=\left[
\begin{array}{cc}
a^{\dagger } & z_{12} \\
z_{21} & z_{22}%
\end{array}%
\right] _{q\times p}.
\end{equation*}%
So,
\begin{equation*}
az=\left[
\begin{array}{cc}
a & 0 \\
0 & 0%
\end{array}%
\right] _{p\times q}\left[
\begin{array}{cc}
a^{\dagger } & z_{12} \\
z_{21} & z_{22}%
\end{array}%
\right] _{q\times p}=\left[
\begin{array}{cc}
p & az_{12} \\
0 & 0%
\end{array}%
\right] _{p\times p}.
\end{equation*}%
Recall thet $p=p^{\ast }$. Thus
\begin{equation*}
\left( az\right) ^{\ast }=\left[
\begin{array}{cc}
p & 0 \\
(az_{12})^{\ast } & 0%
\end{array}%
\right] _{p\times p}.
\end{equation*}%
From $\left( az\right) ^{\ast }=az$ we get $az_{12}=0$. Since $a^{\circ
}=q^{\circ }$, we get $0=qz_{12}=z_{12}$. It follows that $az=\left[
\begin{array}{cc}
p & 0 \\
0 & 0%
\end{array}%
\right] _{p\times p}=p$.
\end{proof}

\begin{corollary}
\label{Corollary1}Let $a\in \mathcal{R}^{\dagger }$ and fix $a^{-\dagger
}\in a\{-\dagger \}$. Then
\begin{equation*}
a\{-\dagger \}=\left\{ a^{-\dagger }+\left( 1-a^{-\dagger }a\right)
waa^{-\dagger }:w\in \mathcal{R}\text{ is arbitrary}\right\} .
\end{equation*}
\end{corollary}

\begin{proof}
Let
\begin{equation*}
\mathcal{S=}\left\{ a^{-\dagger }+\left( 1-a^{-\dagger }a\right)
waa^{-\dagger }:w\in \mathcal{R}\text{ is arbitrary}\right\} .
\end{equation*}%
Let us first show that $\mathcal{S\subseteq }a\{-\dagger \}$. By Lemma \ref%
{Lemma_first}, $aa^{-\dagger }=aa^{\dagger }$ and $aa^{-\dagger }a=a$, and thus
\begin{eqnarray*}
a\left( a^{-\dagger }+\left( 1-a^{-\dagger }a\right) waa^{-\dagger }\right)
&=&aa^{-\dagger }+\left( a-aa^{-\dagger }a\right) waa^{-\dagger } \\
&=&aa^{-\dagger }+\left( a-a\right) waa^{-\dagger } \\
&=&aa^{-\dagger }=aa^{\dagger }.
\end{eqnarray*}%
Also, by Lemma \ref{Lemma_first}, $a^{-\dagger }aa^{-\dagger }=a^{-\dagger }$
and hence

\noindent $\left( a^{-\dagger }+\left( 1-a^{-\dagger }a\right) waa^{-\dagger
}\right) a\left( a^{-\dagger }+\left( 1-a^{-\dagger }a\right) waa^{-\dagger
}\right) =\left( a^{-\dagger }+\left( 1-a^{-\dagger }a\right) waa^{-\dagger
}\right) aa^{-\dagger }=a^{-\dagger }+\left( 1-a^{-\dagger }a\right)
waa^{-\dagger }.$ By Theorem \ref{Theorem1} ((i)$\Leftrightarrow $(ii)) we
have that $\mathcal{S\subseteq }a\{-\dagger \}$.

To prove the reverse inclusion, let $z\in $ $a\{-\dagger \}$. Then there
exists $c\in a\{1\}$ such that $z=caa^{\dagger }$. Note that there exists $%
a^{-}\in a\{1\}$ such that $a^{-\dagger }=a^{-}aa^{\dagger }$, and let $%
p=aa^{-}$ and $q=a^{-}a$. Then $qcp=a^{-}(aca)a^{-}=a^{-}aa^{-}$ and we may
thus write $c=a^{-}aa^{-}+c-$ $qcp$. It follows that
\begin{eqnarray*}
z &=&caa^{\dagger }=\left( a^{-}aa^{-}+c-qcp\right) aa^{\dagger } \\
&=&a^{-}(aa^{-}a)a^{\dagger }+caa^{\dagger }-a^{-}ac(aa^{-}a)a^{\dagger } \\
&=&a^{-\dagger }+\left( 1-a^{-}a\right) caa^{\dagger }.
\end{eqnarray*}%
By Lemma \ref{Lemma_first}, $a^{-\dagger }a=a^{-}a$ and $aa^{-\dagger
}=aa^{\dagger }$. Thus, $z=a^{-\dagger }+\left( 1-a^{-\dagger }a\right)
caa^{-\dagger }$ and so $z\in \mathcal{S}$, i.e. $a\{-\dagger \}\subseteq
\mathcal{S}$.
\end{proof}

We say that $a\in \mathcal{R}^{\dagger }$ is a partial isometry if $a^{\ast
}=a^{\dagger }$ (see \cite{MosicDjordjevic}). We now apply Theorem \ref%
{Theorem1} and its corollary to the subset of all partial isometries in $%
\mathcal{R}$.

\begin{proposition}
Let $a\in \mathcal{R}^\dag$ be a partial isometry. Then the system $xax=x$, $%
ax=aa^{\ast }$ has a solution $x=c$ if and only if $c=a^{-}aa^{\ast
}+(1-a^{-}a)waa^{\ast }$ where $a^{-}\in a\{1\}$ and $w\in \mathcal{R}$.
\end{proposition}

\begin{proof}
Suppose first the system $xax=x$, $ax=aa^{\ast }$ has a solution $x=c$.
Since $a^{\ast }=a^{\dagger }$ we have $cac=c$ and $ac=aa^{\dagger }$. By
Theorem \ref{Theorem1}, $c\in a\{-\dagger \}$, and thus by Corollary \ref%
{Corollary1}, $c=a^{-\dagger }+\left( 1-a^{-\dagger }a\right) waa^{-\dagger
} $ for some $w\in \mathcal{R}$ and $a^{-}\in a\{1\}$ where $a^{-\dagger
}=a^{-}aa^{\dagger }$. By Lemma \ref{Lemma_first}, $a^{-\dagger }a=a^{-}a$,
and $aa^{-\dagger }=aa^{\dagger }=aa^{\ast }$. So, $c=a^{-}aa^{\dagger
}+\left( 1-a^{-}a\right) waa^{\ast }$.

Conversely, let $a^{-}\in a\{1\}$, $w\in \mathcal{R}$ and $c=a^{-}aa^{\ast
}+(1-a^{-}a)waa^{\ast }$. Then
\begin{equation*}
ac=(aa^{-}a)a^{\ast }+(a-aa^{-}a)waa^{\ast }=aa^{\ast }
\end{equation*}%
and
\begin{eqnarray*}
cac &=&\left( a^{-}aa^{\ast }+(1-a^{-}a)waa^{\ast }\right) \left(
aa^{-}aa^{\ast }+\left( a-aa^{-}a\right) waa^{\ast }\right) \\
&=&\left( a^{-}aa^{\ast }+(1-a^{-}a)waa^{\ast }\right) aa^{\ast } \\
&=&a^{-}aa^{\dagger }aa^{\dagger }+(1-a^{-}a)waa^{\dagger }aa^{\dagger } \\
&=&a^{-}aa^{\dagger }+(1-a^{-}a)waa^{\dagger } \\
&=&a^{-}aa^{\ast }+(1-a^{-}a)waa^{\ast }=c,
\end{eqnarray*}%
and hence $x=c$ is a solution of the system $xax=x$, $ax=aa^{\ast }$.
\end{proof}

Several necessary and sufficient conditions for $x=a^{-\dagger }$ are next developed.

\begin{lemma} Let $a\in \mathcal{R}^{\dagger }$ and $a^{-}\in $ $a\{1\}$.
Then the following statements are equivalent:

\begin{itemize}
\item[(i)] $x=a^{-\dagger }$;

\item[(ii)] $ax=aa^\dag$ and $x=a^-ax$;

\item[(iii)] $a^*ax=a^*$ and $x=a^-ax$;

\item[(iv)] $xa=a^-a$ and $x=xaa^\dag$;

\item[(v)] $xaa^-=a^-aa^-$ and $x=xaa^\dag$;

\item[(vi)] $axa=a$ and $a^-axaa^\dag=x$;

\item[(vii)] $xax=x$ and $a^*ax=a^*$.
\end{itemize}
\end{lemma}

\begin{proof} (i) $\Rightarrow$ (ii)--(vii): It can be verified by $x=a^-aa^\dag$.

(ii) $\Rightarrow$ (i): The assumptions $ax=aa^\dag$ and $x=a^-ax$ give $x=a^-(ax)=a^-aa^\dag$.

(iii) $\Rightarrow$ (ii): Multiplying $a^*ax=a^*$ by $(a^\dag)^*$ from the left hand side, we obtain $ax=aa^\dag$.

(iv) $\Rightarrow$ (i): From $xa=a^-a$ and $x=xaa^\dag$, we have $x=(xa)a^\dag=a^-aa^\dag$.

The rest follows similarly.
\end{proof}

We can observe that the existence of 1MP-inverses is closely related with
the existence of adequate idempotents.  If $x\in \mathcal{R}$ satisfies $(ax)^*=ax$, then
it is called a $\{4\}$\textit{%
-inverse} of $a$. If $x\in \mathcal{R}$ satisfies the first two and the forth equations in (\ref{def_MP_inverse}), then
it is called a $\{1,2,4\}$\textit{%
-inverse} of $a$. We denote the sets of all $\{4\}$-inverses and $\{1,2,4\}$%
-inverses of $a\in \mathcal{R}$ by $a\{4\}$ and $a\{1,2,4\}$, respectively.

\begin{theorem}\label{te-1MP} Let $a\{4\}\neq\emptyset$. The following statements are
equivalent:
\begin{itemize}
\item[\rm(i)] $a\{-\dag\}\neq\emptyset$;

\item[\rm(ii)] there exist a projection $p\in\mathcal{R}$ and an idempotent $q\in\mathcal{R}$
such that $p\mathcal{R}=a\mathcal{R}$ and $\mathcal{R} q=\mathcal{R} a$.
\end{itemize}
In addition, for arbitrary $a^-\in a\{1\}$, $qa^-p\in
a\{-\dag\}$, that is, $$q\cdot a\{1\}\cdot p\subseteq
a\{-\dag\}.$$
\end{theorem}

\begin{proof} (i) $\Rightarrow$ (ii): For $x=a^-aa^\dag$, let $p=ax$ and $q=xa$. Since $p=aa^\dag$ and $q=a^-a$, we have
$p=p^2=p^*$, $q=q^2$, $p\mathcal{R}=aa^\dag\mathcal{R}=a\mathcal{R}$ and $\mathcal{R} q=\mathcal{R} a^-a=\mathcal{R} a$.

(ii) $\Rightarrow$ (i): The hypothesis $p\mathcal{R}=a\mathcal{R}$ implies $a=pa$ and $p=au$, for some $u\in\mathcal{R}$.
Now, by $a=aua$, $p^*=p=au$ and $a\{4\}\neq\emptyset$, we deduce that $a^\dag$ exists and thus $a\{-\dag\}\neq\emptyset$. 

Notice that $p=au=aa^\dag au=(auaa^\dag)^*=aa^\dag$. Let $x=qa^-p$, for $a^-\in a\{1\}$.
From $\mathcal{R} q=\mathcal{R} a$, we conclude that
$a=aq$. Therefore, $ax=(aq)a^-p=aa^-p=p=aa^\dag$ and $xax=xp=x$, i.e. $x\in a\{-\dag\}$ by Theorem \ref{Theorem1}.
\end{proof}

For $x, y\in a\{-\dag\}$, we prove that $xay\in a\{-\dag\}$.

\begin{theorem}\label{te2-1MP} Let $a\in\mathcal{R}^\dag$. Then
$$a\{-\dag\}\cdot a\cdot a\{-\dag\}\subseteq a\{-\dag\}.$$
\end{theorem}

\begin{proof} Let $x, y\in a\{-\dag\}$ and $x_1=xay$.
Then, by Theorem \ref{Theorem1}, $ax_1=(axa)y=ay=aa^\dag$ and $x_1ax_1=x_1aa^\dag=xay(aa^\dag)=xa(yay)=xay=x_1$, i.e. $x_1\in a\{-\dag\}$.
So, $a\{-\dag\}\cdot a\cdot a\{-\dag\}\subseteq a\{-\dag\}.$
\end{proof}

\section{The 1MP-partial order}

A partial order associated with 1MP-inverses was introduced on $M_{m,n}(%
\mathbb{C})$ in \cite{Hernandez}. Let $A,B\in M_{m,n}(\mathbb{C})$. We say
that $A$ is below $B$ with respect to the relation $\leq ^{-\dagger }$ and
write $A\leq ^{-\dagger }B$ if there exists a 1MP-inverse $A^{-\dagger }$ of
$A$ such that
\begin{equation*}
A^{-\dagger }A=A^{-\dagger }B\quad \text{and\quad }AA^{-\dagger
}=BA^{-\dagger }.
\end{equation*}
It was proved in \cite{Hernandez} that the relation $\leq ^{-\dagger }$ is a
partial order. We now extend this definition to $\ast $-rings with identity.

\begin{definition}
\label{Definition_1MP_relation}Let $a\in \mathcal{R}^{\dagger }$ and $b\in
\mathcal{R}$. We say that $a$ is below $b$ with respect to the 1MP-relation $%
\leq ^{-\dagger }$ and write $a\leq ^{-\dagger }b$ if there exists $%
a^{-\dagger }\in a\{-\dagger \}$ such that $a^{-\dagger }a=a^{-\dagger }b$
and $aa^{-\dagger }=ba^{-\dagger }$.
\end{definition}

We will prove that the 1MP-relation $\leq ^{-\dagger }$ is a partial order
on $\mathcal{R}^{\dagger }$ for any $\ast $-ring $\mathcal{R}$ with
identity. First, let us present and prove two results that extend Theorem
4.2 in \cite{Hernandez}. With the next theorem we describe the set of all
elements $b\in \mathcal{R}$ that are above a given $a\in \mathcal{R}%
^{\dagger }$ with respect to the 1MP-relation.

\begin{theorem}
\label{Theorem2_0}Let $a\in \mathcal{R}^{\dagger }$ and denote $%
p=aa^{\dagger }$ and $q=a^{\dagger }a$. Let $b\in \mathcal{R}$. Then the
following statements are equivalent:

\begin{itemize}
\item[(i)] $a\leq ^{-\dagger }b$;

\item[(ii)]
\begin{equation}
b=\left[
\begin{array}{cc}
a & 0 \\
-b_{4}da & b_{4}%
\end{array}%
\right] _{p\times q}  \label{Form_of_b_1MP}
\end{equation}
where $b_{4}\in (1-p)\mathcal{R}(1-q)$ and $d\in (1-q)\mathcal{R}p$.
\end{itemize}
\end{theorem}

\begin{proof}
Suppose first $a\leq ^{-\dagger }b$ for $a\in \mathcal{R}^{\dagger }$ and $%
b\in \mathcal{R}$. Then there exists $a^{-\dagger }\in a\{-\dagger \}$ such
that $a^{-\dagger }a=a^{-\dagger }b$ and $aa^{-\dagger }=ba^{-\dagger }$. By
(\ref{minus plus set}) there is $d\in (1-q)\mathcal{R}p$ such that
\begin{equation*}
a^{-\dagger }=\left[
\begin{array}{cc}
a^{\dagger } & 0 \\
d & 0%
\end{array}%
\right] _{q\times p}.
\end{equation*}%
Let $b=\left[
\begin{array}{cc}
b_{1} & b_{2} \\
b_{3} & b_{4}%
\end{array}%
\right] _{p\times q}$. Then
\begin{eqnarray*}
a^{-\dagger }b &=&\left[
\begin{array}{cc}
a^{\dagger } & 0 \\
d & 0%
\end{array}%
\right] _{q\times p}\left[
\begin{array}{cc}
b_{1} & b_{2} \\
b_{3} & b_{4}%
\end{array}%
\right] _{p\times q} \\
&=&\left[
\begin{array}{cc}
a^{\dagger }b_{1} & a^{\dagger }b_{2} \\
db_{1} & db_{2}%
\end{array}%
\right] _{q\times q}
\end{eqnarray*}%
which equals
\begin{equation*}
a^{-\dagger }a=\left[
\begin{array}{cc}
a^{\dagger } & 0 \\
d & 0%
\end{array}%
\right] _{q\times p}\left[
\begin{array}{cc}
a & 0 \\
0 & 0%
\end{array}%
\right] _{p\times q}=\left[
\begin{array}{cc}
a^{\dagger }a & 0 \\
da & 0%
\end{array}%
\right] _{q\times q}.
\end{equation*}%
So, $a^{\dagger }a=a^{\dagger }b_{1}$ and $a^{\dagger }b_{2}=0$, and thus $%
a=aa^{\dagger }b_{1}=pb_{1}$ and $0=aa^{\dagger }b_{2}=pb_{2}$. Since $%
b_{1},b_{2}\in p\mathcal{R}$, we get $a=b_{1}$ and $b_{2}=0$. It follows
that
\begin{equation*}
ba^{-\dagger }=\left[
\begin{array}{cc}
a & 0 \\
b_{3} & b_{4}%
\end{array}%
\right] _{p\times q}\left[
\begin{array}{cc}
a^{\dagger } & 0 \\
d & 0%
\end{array}%
\right] _{q\times p}=\left[
\begin{array}{cc}
p & 0 \\
b_{3}a^{\dagger }+b_{4}d & 0%
\end{array}%
\right] _{p\times p}
\end{equation*}%
which equals $aa^{\dagger }=\left[
\begin{array}{cc}
p & 0 \\
0 & 0%
\end{array}%
\right] _{p\times p}$. So, $b_{3}a^{\dagger }+b_{4}d=0$, i.e. $%
b_{3}a^{\dagger }=-b_{4}d$, and thus $b_{3}q=-b_{4}da$. Since $b_{3}\in
\mathcal{R}q$, we have that $b_{3}=-b_{4}da$ and therefore
\begin{equation*}
b=\left[
\begin{array}{cc}
a & 0 \\
-b_{4}da & b_{4}%
\end{array}%
\right] _{p\times q}.
\end{equation*}

To prove the converse implication, let $a\in \mathcal{R}^{\dagger }$ and let
$b\in \mathcal{R}$ be of matrix form (\ref{Form_of_b_1MP}). Define
\begin{equation*}
c=\left[
\begin{array}{cc}
a^{\dagger } & 0 \\
d & 0%
\end{array}%
\right] _{q\times p}
\end{equation*}%
where the element $d$ is as in (\ref{Form_of_b_1MP}). By (\ref{minus plus
set}), $c\in a\{-\dagger \}$. We have
\begin{equation*}
ca=\left[
\begin{array}{cc}
a^{\dagger } & 0 \\
d & 0%
\end{array}%
\right] _{q\times p}\left[
\begin{array}{cc}
a & 0 \\
0 & 0%
\end{array}%
\right] _{p\times q}=\left[
\begin{array}{cc}
q & 0 \\
da & 0%
\end{array}%
\right] _{q\times q}
\end{equation*}%
and%
\begin{equation*}
cb=\left[
\begin{array}{cc}
a^{\dagger } & 0 \\
d & 0%
\end{array}%
\right] _{q\times p}\left[
\begin{array}{cc}
a & 0 \\
-b_{4}da & b_{4}%
\end{array}%
\right] _{p\times q}=\left[
\begin{array}{cc}
q & 0 \\
da & 0%
\end{array}%
\right] _{q\times q}
\end{equation*}%
and thus $ca=cb$. Also, $ac=aa^{\dagger }=p$ and since $d\in \mathcal{R}p$
we get
\begin{eqnarray*}
bc &=&\left[
\begin{array}{cc}
a & 0 \\
-b_{4}da & b_{4}%
\end{array}%
\right] _{p\times q}\left[
\begin{array}{cc}
a^{\dagger } & 0 \\
d & 0%
\end{array}%
\right] _{q\times p}=\left[
\begin{array}{cc}
p & 0 \\
-b_{4}dp+b_{4}d & 0%
\end{array}%
\right] _{p\times p} \\
&=&\left[
\begin{array}{cc}
p & 0 \\
-b_{4}d+b_{4}d & 0%
\end{array}%
\right] _{p\times p}=\left[
\begin{array}{cc}
p & 0 \\
0 & 0%
\end{array}%
\right] _{p\times p}=p.
\end{eqnarray*}%
So, $ac=bc$ and therefore $a\leq ^{-\dagger }b$.
\end{proof}

For $a,b\in \mathcal{R}^{\dagger }$ where $a\leq ^{-\dagger }b$ the next
result describes the set $b\{-\dagger \}$ of all 1MP-inverses of $b$.

\begin{theorem}
\label{Theorem2}Let $a,b\in \mathcal{R}^{\dagger }$, suppose $a\leq
^{-\dagger }b$, and represent $b$ in the matrix form (\ref{Form_of_b_1MP}).
Then $x\in b\{-\dagger \}$ if and only if there exist $x_{3}\in (1-q)%
\mathcal{R}p$ and $x_{4}\in (1-q)\mathcal{R}(1-p)$ such that
\begin{equation*}
x=\left[
\begin{array}{cc}
a^{\dagger } & 0 \\
x_{3} & x_{4}%
\end{array}%
\right] _{q\times p}
\end{equation*}%
where (for $b_{4}$ and $d$ from (\ref{Form_of_b_1MP})) $b_{4}x_{3}=b_{4}d$
and $x_{4}\in b_{4}\{-\dagger \}$.
\end{theorem}

\begin{proof}
Suppose first $x\in b\{-\dagger \}$ and let
\begin{equation*}
x=\left[
\begin{array}{cc}
x_{1} & x_{2} \\
x_{3} & x_{4}%
\end{array}%
\right] _{q\times p}
\end{equation*}%
for $p=aa^{\dagger }$ and $q=a^{\dagger }a$. By Theorem \ref{Theorem1}, $%
b\{-\dagger \}=b\{1,2,3\}$ and so $bxb=b$, $xbx=x$, and $\left( bx\right) ^{\ast }=bx$. We have
\begin{eqnarray*}
b &=&bxb=\left[
\begin{array}{cc}
a & 0 \\
-b_{4}da & b_{4}%
\end{array}%
\right] _{p\times q}\left[
\begin{array}{cc}
x_{1} & x_{2} \\
x_{3} & x_{4}%
\end{array}%
\right] _{q\times p}\left[
\begin{array}{cc}
a & 0 \\
-b_{4}da & b_{4}%
\end{array}%
\right] _{p\times q} \\
&=&\left[
\begin{array}{cc}
ax_{1}a-ax_{2}b_{4}da & ax_{2}b_{4} \\
\left( -b_{4}dax_{1}+b_{4}x_{3}\right) a-\left(
-b_{4}dax_{2}+b_{4}x_{4}\right) b_{4}da & \left(
-b_{4}dax_{2}+b_{4}x_{4}\right) b_{4}%
\end{array}%
\right] _{p\times q}
\end{eqnarray*}%
and therefore $ax_{2}b_{4}=0$. Also, $ax_{1}a-ax_{2}b_{4}da=a$ and thus $%
ax_{1}a=a$. Hence, $a^{\dagger }ax_{1}aa^{\dagger }=a^{\dagger }aa^{\dagger
}=a^{\dagger }$ and so $qx_{1}p=a^{\dagger }$. Since $x_{1}\in q\mathcal{R}p$%
, it follows that $x_{1}=a^{\dagger }$. The element in the second row and
second column yields $-b_{4}dax_{2}b_{4}+b_{4}x_{4}b_{4}=b_{4}$ and since $%
ax_{2}b_{4}=0$, we have
\begin{equation}
b_{4}x_{4}b_{4}=b_{4}.  \label{eq1}
\end{equation}%
Also, $%
-b_{4}dax_{1}a+b_{4}x_{3}a+b_{4}dax_{2}b_{4}da-b_{4}x_{4}b_{4}da=-b_{4}da$.
Since $ax_{1}a=a$, $ax_{2}b_{4}=0$, and $b_{4}x_{4}b_{4}=b_{4}$, we get $%
b_{4}x_{3}a=b_{4}da$ and thus $b_{4}x_{3}p=b_{4}dp$. From $x_{3},d\in
\mathcal{R}p$ we obtain
\begin{equation*}
b_{4}x_{3}=b_{4}d.
\end{equation*}%
It follows that
\begin{eqnarray*}
bx &=&\left[
\begin{array}{cc}
a & 0 \\
-b_{4}da & b_{4}%
\end{array}%
\right] _{p\times q}\left[
\begin{array}{cc}
a^{\dagger } & x_{2} \\
x_{3} & x_{4}%
\end{array}%
\right] _{q\times p} \\
&=&\left[
\begin{array}{cc}
p & ax_{2} \\
-b_{4}dp+b_{4}x_{3} & -b_{4}dax_{2}+b_{4}x_{4}%
\end{array}%
\right] _{p\times p} \\
&=&\left[
\begin{array}{cc}
p & ax_{2} \\
0 & -b_{4}dax_{2}+b_{4}x_{4}%
\end{array}%
\right] _{p\times p}.
\end{eqnarray*}%
Note that then
\begin{equation*}
\left( bx\right) ^{\ast }=\left[
\begin{array}{cc}
p & 0 \\
(ax_{2})^{\ast } & (-b_{4}dax_{2}+b_{4}x_{4})^{\ast }%
\end{array}%
\right] _{p\times p}
\end{equation*}%
and thus $ax_{2}=0$ which implies $qx_{2}=0$. Since $x_{2}\in q\mathcal{R}$,
we have $x_{2}=0$ and therefore
\begin{equation*}
x=\left[
\begin{array}{cc}
a^{\dagger } & 0 \\
x_{3} & x_{4}%
\end{array}%
\right] _{q\times p}
\end{equation*}%
and $-b_{4}dax_{2}+b_{4}x_{4}=b_{4}x_{4}$. So, $\left( bx\right) ^{\ast
}=bx$ implies
\begin{equation}
\left( b_{4}x_{4}\right) ^{\ast }=b_{4}x_{4}.  \label{eq2}
\end{equation}%
Finally,
\begin{equation*}
x=xbx=\left[
\begin{array}{cc}
a^{\dagger } & 0 \\
x_{3} & x_{4}%
\end{array}%
\right] _{q\times p}\left[
\begin{array}{cc}
a & 0 \\
-b_{4}da & b_{4}%
\end{array}%
\right] _{p\times q}\left[
\begin{array}{cc}
a^{\dagger } & 0 \\
x_{3} & x_{4}%
\end{array}%
\right] _{q\times p}=\left[
\begin{array}{cc}
a^{\dagger } & 0 \\
\ast & x_{4}b_{4}x_{4}%
\end{array}%
\right] _{q\times p}
\end{equation*}%
and thus
\begin{equation}
x_{4}=x_{4}b_{4}x_{4}.  \label{eq3}
\end{equation}%
By equations (\ref{eq1})-(\ref{eq3}) and Theorem \ref{Theorem1} we establish
that $x_{4}\in b_{4}\{-\dagger \}$.

Conversely, let
\begin{equation*}
x=\left[
\begin{array}{cc}
a^{\dagger } & 0 \\
x_{3} & x_{4}%
\end{array}%
\right] _{q\times p}
\end{equation*}%
with $b_{4}x_{3}=b_{4}d$ and $x_{4}\in b_{4}\{-\dagger \}$. Since $%
b_{4}x_{4}b_{4}=b_{4}$, $b_{4}x_{3}=b_{4}d$, and $d\in \mathcal{R}p$, we
obtain
\begin{eqnarray*}
bxb &=&\left[
\begin{array}{cc}
a & 0 \\
-b_{4}da & b_{4}%
\end{array}%
\right] _{p\times q}\left[
\begin{array}{cc}
a^{\dagger } & 0 \\
x_{3} & x_{4}%
\end{array}%
\right] _{q\times p}\left[
\begin{array}{cc}
a & 0 \\
-b_{4}da & b_{4}%
\end{array}%
\right] _{p\times q} \\
&=&\left[
\begin{array}{cc}
a & 0 \\
\left( -b_{4}dp+b_{4}x_{3}\right) a-b_{4}x_{4}b_{4}da & b_{4}x_{4}b_{4}%
\end{array}%
\right] _{p\times q} \\
&=&\left[
\begin{array}{cc}
a & 0 \\
-b_{4}da & b_{4}%
\end{array}%
\right] _{p\times q}=b.
\end{eqnarray*}%
From $x_{4}b_{4}x_{4}=x_{4}$, $b_{4}x_{3}=b_{4}d$, and $x_{3},d\in \mathcal{R%
}p$ we get
\begin{eqnarray*}
xbx &=&\left[
\begin{array}{cc}
a^{\dagger } & 0 \\
x_{3} & x_{4}%
\end{array}%
\right] _{q\times p}\left[
\begin{array}{cc}
a & 0 \\
-b_{4}da & b_{4}%
\end{array}%
\right] _{p\times q}\left[
\begin{array}{cc}
a^{\dagger } & 0 \\
x_{3} & x_{4}%
\end{array}%
\right] _{q\times p} \\
&=&\left[
\begin{array}{cc}
a^{\dagger } & 0 \\
x_{3}p-x_{4}b_{4}dp+x_{4}b_{4}x_{3} & x_{4}b_{4}x_{4}%
\end{array}%
\right] _{q\times p}=\left[
\begin{array}{cc}
a^{\dagger } & 0 \\
x_{3} & x_{4}%
\end{array}%
\right] _{q\times p}=x
\end{eqnarray*}%
and
\begin{eqnarray*}
bx &=&\left[
\begin{array}{cc}
a & 0 \\
-b_{4}da & b_{4}%
\end{array}%
\right] _{p\times q}\left[
\begin{array}{cc}
a^{\dagger } & 0 \\
x_{3} & x_{4}%
\end{array}%
\right] _{q\times p} \\
&=&\left[
\begin{array}{cc}
p & 0 \\
-b_{4}dp+b_{4}x_{3} & b_{4}x_{4}%
\end{array}%
\right] _{p\times p}=\left[
\begin{array}{cc}
p & 0 \\
0 & b_{4}x_{4}%
\end{array}%
\right] _{p\times p}.
\end{eqnarray*}%
Since $p^{\ast }=p$ and $\left( b_{4}x_{4}\right) ^{\ast }=b_{4}x_{4}$, we
establish that $bx=\left( bx\right) ^{\ast }$. Thus, $x\in b\{1,2,3\}$ and
therefore by Theorem \ref{Theorem1}, $x\in b\{-\dagger \}$.
\end{proof}

Let us now prove that the relation $\leq ^{-\dagger }$ is a partial order on
$\mathcal{R}^{\dagger }$.

\begin{theorem}\label{partial_original}
Let $\mathcal{R}$ be a $\ast $-ring with identity. The 1MP-relation $\leq
^{-\dagger }$ is a partial order on $\mathcal{R}^{\dagger }$.
\end{theorem}

\begin{proof}
The set $a\{-\dagger \}$ is nonempty for every $a\in \mathcal{R}^{\dagger}$ which implies that $\leq ^{-\dagger }$ is reflexive. To prove that it is
antisymmetric, let for $a,b\in \mathcal{R}^{\dagger }$, $a\leq ^{-\dagger }b$
and $b\leq ^{-\dagger }a$. Thus, there exists $b^{-\dagger }\in b\{-\dagger
\}$ such that $b^{-\dagger }b=b^{-\dagger }a$ and $bb^{-\dagger
}=ab^{-\dagger }$. Let $p=aa^{\dagger }$ and $q=a^{\dagger }a$. Since $a\leq
^{-\dagger }b$, it follows by Theorems \ref{Theorem2_0} and \ref{Theorem2}
that
\begin{equation}
b=\left[
\begin{array}{cc}
a & 0 \\
-b_{4}da & b_{4}%
\end{array}%
\right] _{p\times q}\quad \text{and\quad }b^{-\dagger }=\left[
\begin{array}{cc}
a^{\dagger } & 0 \\
x_{3} & x_{4}%
\end{array}%
\right] _{q\times p}  \label{Eq_patrial1}
\end{equation}%
where $b_{4}\in (1-p)\mathcal{R}(1-q)$, $d,x_{3}\in (1-q)\mathcal{R}p$, $%
x_{4}\in (1-q)\mathcal{R}(1-p)$, $b_{4}x_{3}=b_{4}d$, and $x_{4}\in
b_{4}\{-\dagger \}$. We have
\begin{equation*}
b^{-\dagger }b=\left[
\begin{array}{cc}
q & 0 \\
x_{3}a-x_{4}b_{4}da & x_{4}b_{4}%
\end{array}%
\right] _{q\times q}
\end{equation*}%
and
\begin{equation*}
b^{-\dagger }a=\left[
\begin{array}{cc}
a^{\dagger } & 0 \\
x_{3} & x_{4}%
\end{array}%
\right] _{q\times p}\left[
\begin{array}{cc}
a & 0 \\
0 & 0%
\end{array}%
\right] _{p\times q}=\left[
\begin{array}{cc}
q & 0 \\
x_{3}a & 0%
\end{array}%
\right] _{q\times q}.
\end{equation*}%
It follows that $x_{4}b_{4}=0$ and thus $b_{4}x_{4}b_{4}=0$. Since $x_{4}\in
b_{4}\{-\dagger \}=b_{4}\{1,2,3\}$, we obtain $b_{4}=$ $b_{4}x_{4}b_{4}=0$ and
therefore
\begin{equation*}
b=\left[
\begin{array}{cc}
a & 0 \\
0 & 0%
\end{array}%
\right] _{p\times q}=a.
\end{equation*}

Let us now prove that $\leq ^{-\dagger }$ is transitive. Let $a\leq
^{-\dagger }b$ and $b\leq ^{-\dagger }c$ for $a,b\in $ $\mathcal{R}^{\dagger
}$ and $c\in \mathcal{R}$. There exists $b^{-\dagger }\in b\{-\dagger \}$
such that $b^{-\dagger }b=b^{-\dagger }c$ and $bb^{-\dagger }=cb^{-\dagger }$%
. Let again $p=aa^{\dagger }$ and $q=a^{\dagger }a$. Since $a\leq ^{-\dagger
}b$, it follows as before by Theorems \ref{Theorem2_0} and \ref{Theorem2}
that $b$ and $b^{-\dagger }$ are of the form (\ref{Eq_patrial1}). Let
\begin{equation*}
c=\left[
\begin{array}{cc}
c_{1} & c_{2} \\
c_{3} & c_{4}%
\end{array}%
\right] _{p\times q}.
\end{equation*}%
From $b^{-\dagger }b=b^{-\dagger }c$ we have
\begin{equation*}
\left[
\begin{array}{cc}
q & 0 \\
x_{3}a-x_{4}b_{4}da & x_{4}b_{4}%
\end{array}%
\right] _{q\times q}=\left[
\begin{array}{cc}
a^{\dagger }c_{1} & a^{\dagger }c_{2} \\
x_{3}c_{1}+x_{4}c_{3} & x_{3}c_{2}+x_{4}c_{4}%
\end{array}%
\right] _{q\times q}
\end{equation*}%
and therefore $a^{\dagger }c_{2}=0$. Thus, $pc_{2}=0$ and so $c_{2}=0$ since
$c_{2}\in p\mathcal{R}$. Also, $q=a^{\dagger }c_{1}$ and hence $aa^{\dagger
}a=aa^{\dagger }c_{1}$, i.e. $a=pc_{1}=c_{1}$. By $bb^{-\dagger
}=cb^{-\dagger }$ we thus obtain
\begin{equation*}
\left[
\begin{array}{cc}
p & 0 \\
-b_{4}dp+b_{4}x_{3} & b_{4}x_{4}%
\end{array}%
\right] _{p\times p}=\left[
\begin{array}{cc}
p & 0 \\
c_{3}a^{\dagger }+c_{4}x_{3} & c_{4}x_{4}%
\end{array}%
\right] _{p\times p}.
\end{equation*}%
So, $-b_{4}dp+b_{4}x_{3}=c_{3}a^{\dagger }+c_{4}x_{3}$. Since $d\in \mathcal{%
R}p$ and $b_{4}x_{3}=b_{4}d$, we have $c_{3}a^{\dagger }=-c_{4}x_{3}$. From
$c_{3}\in \mathcal{R}q$, we establish that $c_{3}=-c_{4}x_{3}a$. It follows
that
\begin{equation*}
c=\left[
\begin{array}{cc}
a & 0 \\
-c_{4}x_{3}a & c_{4}%
\end{array}%
\right] _{p\times q}
\end{equation*}%
and therefore by Theorem \ref{Theorem2_0}, $a\leq ^{-\dagger }c$.
\end{proof}

We give some equivalent conditions for $a\leq ^{-\dagger }b$ to be satisfied.

\begin{lemma}\label{le2-1MP} Let $a\in \mathcal{R}^{\dagger }$ and $b\in \mathcal{R}$.
Then the following statements are equivalent:

\begin{itemize}
\item[(i)] $a\leq ^{-\dagger }b$;

\item[(ii)] there exist $x,y\in a\{-\dagger \}$ such that $xa=xb$
and $ay=by$;

\item[(iii)] there exists $a^{-}\in $ $a\{1\}$ such that $aa^\dag b=a=ba^-a$.
\end{itemize}
\end{lemma}

\begin{proof} (i) $\Rightarrow$ (ii): It is clear.

(ii) $\Rightarrow$ (iii): Since $x=zaa^\dag$ and $y=a^-aa^\dag$, for some $z, a^{-}\in $ $a\{1\}$, we obtain
$a=a(xa)=axb=(aza)a^\dag b=aa^\dag b$ and $a=(ay)a=bya=ba^-(aa^\dag a)=ba^-a$.

(iii) $\Rightarrow$ (i): Using $aa^\dag b=a=ba^-a$, for $x=a^-aa^\dag$, we have
$xa=a^-aa^\dag a=a^-a=a^-aa^\dag b=xb$ and $ax=aa^-aa^\dag=aa^\dag=ba^-aa^\dag=bx$.
Hence, $a\leq ^{-\dagger }b$.
\end{proof}

\begin{theorem} \label{Th_dual1} Let $a\{4\}\neq\emptyset$. The following statements are
equivalent:
\begin{itemize}
\item[\rm(i)] $a\leq ^{-\dagger }b$;

\item[\rm(ii)] there exist a projection $p\in\mathcal{R}$ and an idempotent $q\in\mathcal{R}$
such that $p\mathcal{R}=a\mathcal{R}$, $\mathcal{R} q=\mathcal{R} a$ and $pb=a=bq$.
\end{itemize}
\end{theorem}

\begin{proof} (i) $\Rightarrow$ (ii): Since $ax=bx$ and $xa=xb$ for some $x=a^-aa^\dag$, set $p=ax$ and $q=xa$.
The rest is clear by the proof of Theorem \ref{te-1MP} and $pb=axb=axa=a=bxa=bq$.

(ii) $\Rightarrow$ (i): By Theorem \ref{te-1MP} and $pb=a=bq$, we have $x=qa^-p\in a\{-\dag\}$ and
$ax=(aq)a^-p=aa^-p=bqa^-p=bx$ and $xa=qa^-(pa)=qa^-a=qa^-pb=xb$.
\end{proof}

\begin{theorem} Let $a,b\in\mathcal{R}^\dag$. If $a\leq ^{-\dagger }b$, then
$$b\{-\dag\}\cdot a\cdot b\{-\dag\}\subseteq a\{-\dag\}.$$
\end{theorem}

\begin{proof} Assume that $ax=bx$ and $xa=xb$ for some $x\in a\{-\dag\}$. Let $y, z\in b\{-\dag\}$ and $x_1=zay$.
Then $ay=a(xa)y=axby=(byax)^*=(bybx)^*=(bx)^*=(ax)^*=ax=aa^\dag$ and similarly $az=aa^\dag$, which give $ax_1=azay=aa^\dag aa^\dag=aa^\dag$.
Since $x_1ax_1=zayazay=zaa^\dag aa^\dag ay=zay=x_1$, we deduce that $x_1\in b\{-\dag\}$.
\end{proof}

One of the best known relations on $\mathcal{R}$ that may be defined with generalized inverses is the \textit{minus order} \cite{Hartwig}. We say that $a\in \mathcal{R}^{(1)}$ is below $b\in \mathcal{R}$ with respect to the minus order and write $a\leq ^{-}b$ when there exists $a^{-}\in a\{1\}$ such that
\begin{equation}
a^{-}a=a^{-}b\quad \text{and\quad }aa^{-}=ba^{-}.  \label{def_minus_original}
\end{equation}%
It is known (see  \cite{Hartwig}) that the relation $\leq ^{-}$ is a partial order on $\mathcal{R}^{(1)}$ for any ring $\mathcal{R}$.  

Clearly, since $a\{-\dag\}\subseteq a\{1\}$, $a\leq ^{-\dagger }b$ implies $a\leq^-b$. We study additional conditions for the converse to hold.

\begin{theorem} \label{Th_dual3} Let $a\in \mathcal{R}^{\dagger }$ and $b\in \mathcal{R}$.
Then the following statements are equivalent:

\begin{itemize}
\item[(i)] $a\leq ^{-\dagger }b$;

\item[(ii)] $a\leq^-b$ and $a^\dag b=a^\dag a$;

\item[(iii)] $aa^-=ba^-$ and $a^\dag b=a^\dag a$, for some $a^{-}\in $ $a\{1\}$.
\end{itemize}
\end{theorem}

\begin{proof} (i) $\Rightarrow$ (ii): Applying (\ref{Matrix form_basic}) and (\ref{Form_of_b_1MP}), notice that
$$a^\dag b=\left[
\begin{array}{cc}
a^{\dagger } & 0 \\
0 & 0%
\end{array}%
\right] _{q\times p}\left[
\begin{array}{cc}
a & 0 \\
-b_{4}da & b_{4}%
\end{array}%
\right] _{p\times q}=\left[
\begin{array}{cc}
q & 0 \\
0 & 0%
\end{array}%
\right] _{q\times q}=a^\dag a.$$

(ii) $\Rightarrow$ (iii): It is evident.

(iii) $\Rightarrow$ (i): Suppose that $aa^-=ba^-$ and $a^\dag b=a^\dag a$, for some $a^{-}\in $ $a\{1\}$.
Then $aa^\dag b=a=ba^-a$ implies $a\leq ^{-\dagger }b$, by Lemma \ref{le2-1MP}.
\end{proof}

\section{The dual case}

A dual version of partial order $\leq ^{-\dagger }$ was introduced in \cite%
{Hernandez} on $M_{m,n}(\mathbb{C})$. Let $A,B\in M_{m,n}(\mathbb{C})$. We
say that $A$ is below $B$ with respect to the relation $\leq ^{\dagger -}$
and write $A\leq ^{\dagger -}B$ if there exists a MP1-inverse $A^{\dagger -}$
of $A$ such that
\begin{equation*}
A^{\dagger -}A=A^{\dagger -}B\quad \text{and\quad }AA^{\dagger
-}=BA^{\dagger -}.
\end{equation*}%
We now extend the concept of MP1-inverses and the MP1-relation $\leq ^{\dagger
-} $ to the setting of $\ast $-rings. Recall that $\mathcal{R}$ denotes a $%
\ast $-ring with identity.

\begin{definition}
Let $a\in \mathcal{R}^{\dagger }$. For each $a^{-}\in $ $a\{1\}$, the
element
\begin{equation*}
a^{\dagger -}=a^{\dagger }aa^{-}
\end{equation*}%
is called a MP1-inverse of $a$. The set of all MP1-inverses of $a$ is
denoted by $a\{\dagger -\}$.
\end{definition}

A MP1-inverse of $a\in \mathcal{R}^{\dagger }$ always exists however in
general it is not necessarily unique.

\begin{definition}
\label{Definition_MP1_relation}Let $a\in \mathcal{R}^{\dagger }$ and $b\in
\mathcal{R}$. We say that $a$ is below $b$ with respect to the MP1-relation $%
\leq ^{\dagger -}$ and write $a\leq ^{\dagger -}b$ if there exists $%
a^{\dagger -}\in a\{\dagger -\}$ such that $a^{\dagger -}a=a^{\dagger -}b$
and $aa^{\dagger -}=ba^{\dagger -}$
\end{definition}

Consider now a new $\ast $-ring $\mathcal{R}_{L}=\left( \mathcal{R},+,\cdot
_{L},\ast \right) $ where for $a,b\in \mathcal{R}$,
\begin{equation}
a\cdot _{L}b:=ba.  \label{Multiplication_dual}
\end{equation}%
Let $a,c\in \mathcal{R}$. Then $aca=a$ if and only if $a\cdot _{L}c\cdot
_{L}a=aca=a$. So, $c=a^{-}$ is an inner generalized inverse of $a$ in $%
\mathcal{R}$ if and only if $c=a^{-}$ is an inner generalized inverse of $a$
in $\mathcal{R}_{L}$. Similarly we observe that for $a\in \mathcal{R}%
^{\dagger }$, its unique Moore-Penrose inverse $a^{\dagger }$ in $\mathcal{R}
$ is also its unique Moore-Penrose inverse in $\mathcal{R}_{L}$. Moreover,
since $\left( a\cdot _{L}c\right) ^{\ast }=a\cdot _{L}c$ if and only if $%
\left( ca\right) ^{\ast }=ca$, we establish that the set $a\{1,2,3\}$ of all
$\{1,2,3\}$-inverses of $a$ in $\mathcal{R}_{L}$ is the set $a\{1,2,4\}$ of
all $\{1,2,4\}$-inverses of $a$ in $\mathcal{R}$. From
\begin{equation*}
a^{-}\cdot _{L}a\cdot _{L}a^{\dagger }=a^{\dagger }aa^{-}
\end{equation*}%
we have that $a^{-\dagger }$ is a 1MP-invrese of $a$ in $\mathcal{R}_{L}$ if
and only if it is a MP1-inverse of $a$ in $\mathcal{R}$, i.e. $z\in
a\{-\dagger \}$ in $\mathcal{R}_{L}$ if and only if $z\in a\{\dagger -\}$ in
$\mathcal{R}$. The following theorem and its corollary are thus a direct
corollary (\ref{Multiplication_dual}), Theorem \ref{Theorem1} and Corollary %
\ref{Corollary1}.

\begin{theorem}
\label{Theorem 3}Let $a\in \mathcal{R}^{\dagger }$. Then the following
statements are equivalent:

\begin{itemize}
\item[(i)] $z\in a\{\dagger -\}$;

\item[(ii)] $z\in \mathcal{R}$ is a solution of the system $xax=x$ and $%
xa=a^{\dagger }a$;

\item[(iii)] $z\in a\{1,2,4\}$.
\end{itemize}
\end{theorem}

\begin{corollary}
\bigskip Let $a\in \mathcal{R}^{\dagger }$ and fix $a^{\dagger -}\in
a\{\dagger -\}$. Then
\begin{equation*}
a\{\dagger -\}=\left\{ a^{\dagger -}+a^{\dagger -}aw\left( 1-aa^{\dagger
-}\right) :w\in \mathcal{R}\text{ is arbitrary}\right\} .
\end{equation*}
\end{corollary}

\begin{remark}
Let $a\in \mathcal{R}^{\dagger }$. By Theorems \ref{Theorem1} and \ref%
{Theorem 3} we have $a\{-\dagger \}=a\{1,2,3\}$ and $a\{\dagger
-\}=a\{1,2,4\}$. It follows that
\begin{equation*}
a\{-\dagger \}\cap a\{\dagger -\}=\mathcal{R}^{\dagger }.
\end{equation*}
\end{remark}

It is easy to see (compare Definitions \ref{Definition_1MP_relation} and \ref%
{Definition_MP1_relation}) that for every $a\in \mathcal{R}^{\dagger }$
and $b\in \mathcal{R}$ we have
\begin{equation*}
a\leq ^{\dagger -}b\quad \text{if and only if}\quad a\leq _{L}^{-\dagger }b
\end{equation*}%
where $\leq _{L}^{-\dagger }$ is the 1MP-relation in $(\mathcal{R}%
,+,\cdot _{L},\ast )$. By Theorem \ref{partial_original} we may immediately establish that MP1-relation $\leq
^{\dagger -}$ introduced with Definition \ref{Definition_MP1_relation} is a
partial order on $\mathcal{R}^{\dagger }$. We conclude this section with two
results that are direct corollaries of (\ref{Multiplication_dual}) and
Theorems \ref{Theorem2_0} and \ref{Theorem2}, respectively.

\begin{theorem}
Let $a\in \mathcal{R}^{\dagger }$ and denote $p=a^{\dagger }a$ and $%
q=aa^{\dagger }$. Let $b\in \mathcal{R}$. Then the following statements are
equivalent:

\begin{itemize}
\item[(i)] $a\leq ^{\dagger -}b$;

\item[(ii)]
\begin{equation}
b=\left[
\begin{array}{cc}
a & -adb_{4} \\
0 & b_{4}%
\end{array}%
\right] _{p\times q}  \label{Form_of_b_MP1}
\end{equation}%
where $b_{4}\in (1-p)\mathcal{R}(1-q)$ and $d\in q\mathcal{R}(1\mathcal{-}p)$%
.
\end{itemize}
\end{theorem}

\begin{theorem}
Let $a,b\in \mathcal{R}^{\dagger }$, suppose $a\leq ^{\dagger -}b$, and
represent $b$ in the matrix form (\ref{Form_of_b_MP1}). Then $x\in
b\{\dagger -\}$ if and only if there exist $x_{3}\in q\mathcal{R}(1-p)$ and $%
x_{4}\in (1-q)\mathcal{R}(1-p)$ such that
\begin{equation*}
x=\left[
\begin{array}{cc}
a^{\dagger } & x_{3} \\
0 & x_{4}%
\end{array}%
\right] _{q\times p}
\end{equation*}%
where $x_{3}b_{4}=db_{4}$ and $x_{4}\in b_{4}\{\dagger -\}$.
\end{theorem}

Let us end this section with a note that we may similarly obtain ``the dual version'' of  Lemma \ref{le2-1MP} and Theorems \ref{Th_dual1}--\ref{Th_dual3}.

\section{The plus partial order in Rickart $\ast $-rings}

A ring $\mathcal{R}$ is called a \textit{Rickart ring} if for every $a\in
\mathcal{R}$ there exist idempotent elements $p,q\in \mathcal{R}$ such that $%
a^{\circ }=p\cdot \mathcal{R}$ and $^{\circ }a=\mathcal{R\cdot }q$. A $\ast $%
-ring $\mathcal{R} $ is a \textit{Rickart }$\ast $\textit{-ring} if the left
annihilator $^{\circ }a$ of any element $a\in \mathcal{R}$ is generated by a
projection $e\in \mathcal{R}$, i.e. $^{\circ }a$ $=\mathcal{%
R\cdot }e$ where $e=e^{\ast }=e^{2}$. The projection $e$ is
unique and every Rickart ring has the (multiplicative) identity $1$ (see
\cite{Baer} or \cite{Kaplansky}). Let us denote
\begin{equation*}
\text{LP}(a)=\left\{ p\in \mathcal{R}:p=p^{2},\text{ }^{\circ }a=\,^{\circ
}p\right\} \quad \text{and\quad RP}(a)=\left\{ q\in \mathcal{R}:q=q^{2},%
\text{ }a^{\circ }=q^{\circ }\right\} .
\end{equation*}%
Note that the sets LP$(a)$ and RP$(a)$ are nonempty in case when $\mathcal{R}
$ is a Rickart ring. If $\mathcal{R}$ is a Rickart $\ast $-ring, there
exists the unique projection in $\text{LP}(a)$. We denote it by
$\text{lp}(a)$. Similarly, let $\text{rp}(a)$ denote the unique projection in $\text{RP}(a)$

Let $\mathcal{H},\mathcal{K}$ be Hilbert spaces and let $\mathcal{B(H},%
\mathcal{K)}$ be the set of all bounded linear operators from $\mathcal{H}$
to $\mathcal{K}$. When $\mathcal{H}=\mathcal{K}$, we write $\mathcal{B(H)}$
insead of $\mathcal{B(H},\mathcal{H)}$. Note that $\mathcal{B(H)}$ is an
example of a Rickart $\ast $-ring. Let $\mathrm{Ker}\,A$, $\func{Im}\,A$, and $\overline{\func{Im}\,A}$ denote the kernel, the range, and the closure of the range of $A\in
\mathcal{B(H},\mathcal{K)}$, respectively. A new
relation was introduced in \cite{Arias} on $\mathcal{B(H},\mathcal{K)}$. The
definition follows.

\begin{definition}
\label{Definition_plus_original}Let $A,B\in \mathcal{B(H},\mathcal{K)}$. We
say that $A$ is below $B$ with respect to the relation $\leq ^{+}$ and write
$A\leq ^{+}B$ if $\func{Im}\,A\subseteq \func{Im}\,B$, $\func{Im}\,A^{\ast } \subseteq
\func{Im}\,B^{\ast } $, and there are idempotent operators $\widetilde{Q}%
\in \mathcal{B(H)}$ and $Q\in \mathcal{B(K)}$ such that $A=\widetilde{Q}BQ$.
In such case, we can without loss of generality assume that $\func{Im}\,\widetilde{Q} =\overline{\func{Im}\,A}$ and $\func{Im}\,Q^{\ast } =
\overline{\func{Im}\,A^{\ast }}$.
\end{definition}

It was proved in \cite{Arias} that $\leq ^{+}$ is a partial order on $%
\mathcal{B(H},\mathcal{K)}$. Authors of \cite{Arias} named the relation $%
\leq ^{+}$ the plus (partial) order. The plus partial order emerges on $%
\mathcal{B(H},\mathcal{K)}$ as a generalization of the minus order and the \textit{diamond order%
} \cite{BaksalaryHauke}.
It is known that for any von Neumann regular ring $\mathcal{R}$ the definition of the minus (partial) order (\ref%
{def_minus_original}) is equivalent to
\begin{equation}
a=pb=bq  \label{minus_def}
\end{equation}%
where $p,q\in \mathcal{R}$ are some idempotents. It turns out (see \cite[%
Corollary 2.1]{Marovt} and \cite{DjordjevicRakicMarovt, DolinarKuzmaMarovt})
that the relation $\leq ^{-}$ defined with (\ref{minus_def}) (i.e. $a\leq
^{-}b$ when $a=pb=bq$ for some $p=p^{2}$ and $q=q^{2}$) is a partial order
when $\mathcal{R}$ is a Rickart ring.

Let $\mathcal{R}$ be a $\ast $-ring. We say that $a$ is below $b$ with
respect to the diamond relation $\leq ^{\diamond }$ and write $a\leq
^{\diamond }b$ when
\begin{equation}
^{\circ }b\subseteq \,^{\circ }a,b^{\circ }\subseteq a^{\circ },\quad \text{%
and\quad }ab^{\ast }a=aa^{\ast }a.  \label{Def_diamond}
\end{equation}%
It turns out (see \cite{Lebtahi, MarovtRakicDjodjevic}) that (at least) when $%
\mathcal{R}$ is $\ast $\textit{-}regular ring with identity, the relation $%
\leq ^{\diamond }$ is a partial order.

For $A,B\in \mathcal{B}(\mathcal{H})$ we have (see \cite{Marovt,
MarovtRakicDjodjevic})
\begin{equation*}
A^{\circ }\subseteq B^{\circ }\text{\quad if and only if\quad }\mathrm{Ker}%
\,A\subseteq \mathrm{Ker}\,B\text{\quad if and only if\quad }\overline{\func{%
Im}\,B^{\ast }}\subseteq \overline{\func{Im}\,A^{\ast }}
\end{equation*}%
and%
\begin{equation*}
{^{\circ }}A\subseteq {^{\circ }}B\text{\quad if and only if\quad }\overline{%
\func{Im}\,B}\subseteq \overline{\func{Im}\,A}.
\end{equation*}%
Motivated by Definition \ref{Definition_plus_original} this observation
leads us to the following definition.

\begin{definition}
\label{Definition_plus_ring}Let $\mathcal{R}$ be a Rickart $\ast $-ring and
let $a,b\in \mathcal{R}$. We say that $a$ is below $b$ with respect to the
plus order $\leq ^{+}$ and write $a\leq ^{+}b$ if $^{\circ }b\subseteq
\,^{\circ }a$, $b^{\circ }\subseteq a^{\circ }$, and there exist $q\in $RP$%
(a),\widetilde{q}\in $LP$(a)$ such that $a=\widetilde{q}bq$.
\end{definition}

\begin{theorem}
Let $\mathcal{R}$ be a Rickart $\ast $-ring. The relation $\leq ^{+}$
introduced with Definition \ref{Definition_plus_ring} is a partial order on $%
\mathcal{R}$.
\end{theorem}

\begin{proof}
Since $\mathcal{R}$ be a Rickart $\ast $-ring, there exist $q\in $RP$(a)$
and $\widetilde{q}\in $LP$(a)$. Note that for every idempotent $p$, $%
p(1-p)=(1-p)p=0$. Since $q\in $RP$(a),\widetilde{q}\in $LP$(a)$, we thus
have $\left( 1-\widetilde{q}\right) a=0=a\left( 1-q\right) $ and hence $a=%
\widetilde{q}a=aq$. So, $a=\widetilde{q}aq$ and therefore $\leq ^{+}$ is
reflexive.

To show that $\leq ^{+}$ is antisymmetric, let $a\leq ^{+}b$ and $b\leq
^{+}a $. Then there exist $q\in $RP$(a)$, $q_{1}\in $RP$(b)$, $\widetilde{q}%
\in $LP$(a)$, and $\widetilde{q}_{1}\in $LP$(b)$ such that $a=\widetilde{q}%
bq $ and $b=\widetilde{q}_{1}aq_{1}$. Also, $^{\circ }b=\,^{\circ }a$ and $%
b^{\circ }=a^{\circ }$, and therefore $\left( 1-\widetilde{q}\right)
\widetilde{q}_{1}=0$ and $q_{1}\left(1-q\right) =0$. So, $\widetilde{q}_{1}=%
\widetilde{q}\widetilde{q}_{1}$ and $q_{1}=q_{1}q$. It follows that
\begin{equation*}
a=\widetilde{q}bq=\widetilde{q}\widetilde{q}_{1}aq_{1}q=\widetilde{q}%
_{1}aq_{1}=b.
\end{equation*}

Let us now prove that $\leq ^{+}$ is transitive. Let for $a,b,c\in \mathcal{R%
}$, $a\leq ^{+}b$ and $b\leq ^{+}c$. So, $^{\circ }b\subseteq \,^{\circ }a$,
$b^{\circ }\subseteq a^{\circ }$ , $^{\circ }c\subseteq \,^{\circ }b$, $%
c^{\circ }\subseteq b^{\circ }$ and therefore $^{\circ }c\subseteq \,^{\circ
}a$ and $c^{\circ }\subseteq a^{\circ }$. Also, there exist $q_{1}\in $RP$%
(a) $, $q_{2}\in $RP$(b)$, $\widetilde{q}_{1}\in $LP$(a)$, and $\widetilde{q}%
_{2}\in $LP$(b)$ with $a=\widetilde{q}_{1}bq_{1}$ and $b=\widetilde{q}%
_{2}cq_{2}$. It follows that $a=\widetilde{q}_{1}\widetilde{q}%
_{2}cq_{2}q_{1} $. Denote $q_{3}=q_{2}q_{1}$ and $\widetilde{q}_{3}=%
\widetilde{q}_{1}\widetilde{q}_{2}$. So,
\begin{equation*}
a=\widetilde{q}_{3}cq_{3}.
\end{equation*}%
To conclude the proof let us show that $q_{3}\in $RP$(a)$ and $\widetilde{q}%
_{3}\in $LP$(a)$. From $^{\circ }b\subseteq \,^{\circ }a$ and $b^{\circ
}\subseteq a^{\circ }$ we get $^{\circ }\widetilde{q}_{2}\subseteq \,^{\circ
}\widetilde{q}_{1}$ and $q_{2}^{\circ }\subseteq q_{1}^{\circ }$, and thus $%
\left( 1-\widetilde{q}_{2}\right) \widetilde{q}_{1}=0$ and $q_{1}\left(
1-q_{2}\right) =0$. So, $\widetilde{q}_{1}=\widetilde{q}_{2}\widetilde{q}%
_{1} $ and $q_{1}=q_{1}q_{2}$. It follows
\begin{equation*}
\widetilde{q}_{3}^{2}=\widetilde{q}_{1}\left( \widetilde{q}_{2}\widetilde{q}%
_{1}\right) \widetilde{q}_{2}=\widetilde{q}_{1}^{2}\widetilde{q}_{2}=%
\widetilde{q}_{1}\widetilde{q}_{2}=\widetilde{q}_{3}
\end{equation*}%
and
\begin{equation*}
q_{3}^{2}=q_{2}\left( q_{1}q_{2}\right)
q_{1}=q_{2}q_{1}^{2}=q_{2}q_{1}=q_{3}.
\end{equation*}%
Observe that from $a=\widetilde{q}_{3}cq_{3}$ we obtain $^{\circ }\widetilde{%
q}_{3}\subseteq \,^{\circ }a$ and $q_{3}^{\circ }\subseteq a^{\circ }$.
Since $\widetilde{q}_{3}=\widetilde{q}_{1}\widetilde{q}_{2}$ and $%
q_{3}=q_{2}q_{1}$, we have $^{\circ }a=$ $^{\circ }\widetilde{q}%
_{1}\subseteq \,^{\circ }\widetilde{q}_{3}$ and $a^{\circ }=q_{1}^{\circ
}\subseteq q_{3}^{\circ }$. Thus $^{\circ }\widetilde{q}_{3}=\,^{\circ }a$
and $q_{3}^{\circ }=a^{\circ }$ and hence $q_{3}\in $RP$(a)$ and $\widetilde{%
q}_{3}\in $LP$(a)$. Therefore, $a\leq ^{+}c$.
\end{proof}

For the rest of the paper, let $\mathcal{R}$ be a Rickart $\ast $-ring. Let $%
a,b\in \mathcal{R}$ with $a\leq ^{\diamond }b$. Then $^{\circ }b\subseteq
\,^{\circ }a,b^{\circ }\subseteq a^{\circ }$, and $ab^{\ast }a=aa^{\ast }a$.
It is easy to prove (see \cite[proof of Lemma 5]{MarovtMihelic}) that
\begin{equation*}
ab^{\ast }a=aa^{\ast }a\quad \text{if and only if\quad }a=\text{lp}(a)b\text{%
rp}(a).
\end{equation*}%
Since lp$(a)\in $LP$(a)$ and rp$(a)\in $RP$(a)$, it follows that $a\leq
^{+}b $. Let now $a\leq ^{-}b$, i.e. $a=pb=bq$ for some idempotents $p,q\in
\mathcal{R}$. Then $a=pbq$, $^{\circ }b\subseteq \,^{\circ }a$, and $%
b^{\circ }\subseteq a^{\circ }$. Also, by \cite[Corollary 2.1]{Marovt} we
may without loss of generality assume that $^{\circ }a=\,^{\circ }p$ and $%
a^{\circ }=q^{\circ }$, i.e. $p\in $LP$(a)$ and $q\in $RP$(a)$. It follows
again that $a\leq ^{+}b$. We sum up these observations in the following
proposition.

\begin{proposition}
Let $a,b\in \mathcal{R}$. If $a\leq ^{\diamond }b$, then $a\leq ^{+}b$, and
if $a\leq ^{-}b$, then $a\leq ^{+}b$.
\end{proposition}

We end the paper with a characterization of the plus partial order in
Rickart $\ast $-rings. First, let us present an auxiliary result which was
proved in \cite{DjordjevicRakicMarovt}.

\begin{lemma}
\label{Lemma_final}Let $a\in \mathcal{R}$, $p\in $LP$(a)$, and $q\in $RP$(a)$%
. Then
\begin{equation*}
\text{LP}(a)=\left\{ \left[
\begin{array}{cc}
p & p_{1} \\
0 & 0%
\end{array}%
\right] _{p\times p}:p_{1}\in p\mathcal{R}(1-p)\right\} \text{ and RP}%
(a)=\left\{ \left[
\begin{array}{cc}
q & 0 \\
q_{1} & 0%
\end{array}%
\right] _{q\times q}:q_{1}\in (1-q)\mathcal{R}q\right\} .
\end{equation*}
\end{lemma}

Let $a\in \mathcal{R}$. For the purposes of the following characterization
of the plus partial order we denote $l_{a}=$lp$(a)$ and $r_{a}=$ rp$(a)$.
The next theorem extends \cite[Theorem 3.11]{Arias} to the setting of
Rickart $\ast $-rings.

\begin{theorem}
Let $\mathcal{R}$ be a Rickart $\ast $-ring and let $a,b\in \mathcal{R}$.
The following statements are equivalent:

\begin{itemize}
\item[(i)] $a\leq ^{+}b;$

\item[(ii)] $b=\left[
\begin{array}{cc}
a+y(b_{22}x+w)+zx & yb_{22}+z \\
b_{22}x+w & b_{22}%
\end{array}%
\right] _{l_{a}\times r_{a}}$\newline
where $b_{22}\in (1-l_{a})\mathcal{R}(1-r_{a})$ is arbitrary, and $y\in l_{a}%
\mathcal{R}(1-l_{a})$, $x\in (1-r_{a})\mathcal{R}r_{a}$, $w\in (1-l_{a})%
\mathcal{R}r_{a}$, and $z\in l_{a}\mathcal{R}(1-r_{a})$ are such that $%
^{\circ }b\subseteq \,^{\circ }\left( \left( y+1\right) w\right) $ and $%
b^{\circ }\subseteq \left( z\left( x+1\right) \right) ^{\circ }$.
\end{itemize}

\begin{proof}
Let $a,b\in \mathcal{R}$. Suppose $y\in l_{a}\mathcal{R}(1-l_{a})$, $x\in
(1-r_{a})\mathcal{R}r_{a}$, and
\begin{equation*}
\widetilde{q}=\left[
\begin{array}{cc}
l_{a} & -y \\
0 & 0%
\end{array}%
\right] _{l_{a}\times l_{a}}\quad \text{and\quad }q=\left[
\begin{array}{cc}
r_{a} & 0 \\
-x & 0%
\end{array}%
\right] _{r_{a}\times r_{a}}.
\end{equation*}%
By Lemma \ref{Lemma_final}, $\widetilde{q}\in $LP$(a)$ and $q\in $RP$(a)$.
We have $a=\widetilde{q}bq$ if and only if
\begin{equation}
b=a+\left( 1-\widetilde{q}\right) b\left( 1-q\right) +\left( 1-\widetilde{q}%
\right) bq+\widetilde{q}b\left( 1-q\right) .  \label{eq_final}
\end{equation}%
Let $b=\left[
\begin{array}{cc}
b_{11} & b_{12} \\
b_{21} & b_{22}%
\end{array}%
\right] _{l_{a}\times r_{a}}$ and denote $w=b_{21}-b_{22}x\in (1-l_{a})%
\mathcal{R}r_{a}$ and $z=b_{12}-yb_{22}\in l_{a}\mathcal{R}(1-r_{a})$. Then
\begin{eqnarray*}
\left( 1-\widetilde{q}\right) b\left( 1-q\right) &=&\left( \left[
\begin{array}{cc}
l_{a} & 0 \\
0 & 1-l_{a}%
\end{array}%
\right] _{l_{a}\times l_{a}}-\left[
\begin{array}{cc}
l_{a} & -y \\
0 & 0%
\end{array}%
\right] _{l_{a}\times l_{a}}\right) \left[
\begin{array}{cc}
b_{11} & b_{12} \\
b_{21} & b_{22}%
\end{array}%
\right] _{l_{a}\times r_{a}}\cdot \\
&&\cdot \left( \left[
\begin{array}{cc}
r_{a} & 0 \\
0 & 1-r_{a}%
\end{array}%
\right] _{r_{a}\times r_{a}}-\left[
\begin{array}{cc}
r_{a} & 0 \\
-x & 0%
\end{array}%
\right] _{r_{a}\times r_{a}}\right) \\
&=&\left[
\begin{array}{cc}
yb_{21} & yb_{22} \\
b_{21} & b_{22}%
\end{array}%
\right] _{l_{a}\times r_{a}}\left[
\begin{array}{cc}
0 & 0 \\
x & 1-r_{a}%
\end{array}%
\right] _{r_{a}\times r_{a}}=\left[
\begin{array}{cc}
yb_{22}x & yb_{22} \\
b_{22}x & b_{22}%
\end{array}%
\right] _{l_{a}\times r_{a}},
\end{eqnarray*}

\begin{eqnarray*}
\left( 1-\widetilde{q}\right) bq &=&\left[
\begin{array}{cc}
0 & y \\
0 & 1-l_{a}%
\end{array}%
\right] _{l_{a}\times l_{a}}\left[
\begin{array}{cc}
b_{11} & b_{12} \\
b_{21} & b_{22}%
\end{array}%
\right] _{l_{a}\times r_{a}}\left[
\begin{array}{cc}
r_{a} & 0 \\
-x & 0%
\end{array}%
\right] _{r_{a}\times r_{a}} \\
&=&\left[
\begin{array}{cc}
yb_{21} & yb_{22} \\
b_{21} & b_{22}%
\end{array}%
\right] _{l_{a}\times r_{a}}\left[
\begin{array}{cc}
r_{a} & 0 \\
-x & 0%
\end{array}%
\right] _{r_{a}\times r_{a}} \\
&=&\left[
\begin{array}{cc}
y\left( b_{21}-b_{22}x\right) & 0 \\
b_{21}-b_{22}x & 0%
\end{array}%
\right] _{l_{a}\times r_{a}}=\left[
\begin{array}{cc}
yw & 0 \\
w & 0%
\end{array}%
\right] _{l_{a}\times r_{a}},
\end{eqnarray*}%
and
\begin{eqnarray*}
\widetilde{q}b\left( 1-q\right) &=&\left[
\begin{array}{cc}
l_{a} & -y \\
0 & 0%
\end{array}%
\right] _{l_{a}\times l_{a}}\left[
\begin{array}{cc}
b_{11} & b_{12} \\
b_{21} & b_{22}%
\end{array}%
\right] _{l_{a}\times r_{a}}\left[
\begin{array}{cc}
0 & 0 \\
x & 1-r_{a}%
\end{array}%
\right] _{r_{a}\times r_{a}} \\
&=&\left[
\begin{array}{cc}
b_{11}-yb_{21} & b_{12}-yb_{22} \\
0 & 0%
\end{array}%
\right] _{l_{a}\times r_{a}}\left[
\begin{array}{cc}
0 & 0 \\
x & 1-r_{a}%
\end{array}%
\right] _{r_{a}\times r_{a}} \\
&=&\left[
\begin{array}{cc}
\left( b_{12}-yb_{22}\right) x & b_{12}-yb_{22} \\
0 & 0%
\end{array}%
\right] _{l_{a}\times r_{a}}=\left[
\begin{array}{cc}
zx & z \\
0 & 0%
\end{array}%
\right] _{l_{a}\times r_{a}}.
\end{eqnarray*}%
By (\ref{eq_final}) it follows that $a=\widetilde{q}bq$ if and only if
\begin{eqnarray*}
b &=&\left[
\begin{array}{cc}
a & 0 \\
0 & 0%
\end{array}%
\right] _{l_{a}\times r_{a}}+\left[
\begin{array}{cc}
yb_{22}x & yb_{22} \\
b_{22}x & b_{22}%
\end{array}%
\right] _{l_{a}\times r_{a}}+\left[
\begin{array}{cc}
yw & 0 \\
w & 0%
\end{array}%
\right] _{l_{a}\times r_{a}}+\left[
\begin{array}{cc}
zx & z \\
0 & 0%
\end{array}%
\right] _{l_{a}\times r_{a}} \\
&=&\left[
\begin{array}{cc}
a+y\left( b_{22}x+w\right) +zx & yb_{22}+z \\
b_{22}x+w & b_{22}%
\end{array}%
\right] _{l_{a}\times r_{a}}.
\end{eqnarray*}

From now on we assume that $a=\widetilde{q}bq$. To finish the proof, let us
show that $^{\circ }b\subseteq \,^{\circ }a$ if and only if $^{\circ
}b\subseteq \,^{\circ }\left( \left( y+1\right) w\right) $, and $b^{\circ
}\subseteq a^{\circ }$ if and only if and $b^{\circ }\subseteq \left(
z\left( x+1\right) \right) ^{\circ }$. We will prove only the first
equivalence, the second equivalence can be proved similarly. Note first that
on the one hand
\begin{equation*}
\left( 1-\widetilde{q}\right) bq=\left[
\begin{array}{cc}
yw & 0 \\
w & 0%
\end{array}%
\right] _{l_{a}\times r_{a}}=yw+w=(y+1)w
\end{equation*}%
and on the other hand $\left( 1-\widetilde{q}\right) bq=bq-\widetilde{q}%
bq=bq-a$. Let now $^{\circ }b\subseteq \,^{\circ }a$ and suppose $tb=0$ for
some $t\in \mathcal{R}$. Thus $ta=0$ and therefore
\begin{equation*}
0=t\left( bq-a\right) =t\left( 1-\widetilde{q}\right) bq=t(y+1)w.
\end{equation*}%
So, $^{\circ }b\subseteq \,^{\circ }\left( \left( y+1\right) w\right) $.
Conversely, let $^{\circ }b\subseteq \,^{\circ }\left( \left( y+1\right)
w\right) $ and $tb=0$. Then $0=t(y+1)w$ and thus $0=t\left( bq-a\right) =-ta$%
. So, $ta=0$ and hence $^{\circ }b\subseteq \,^{\circ }a$.
\end{proof}
\end{theorem}

\end{document}